\documentclass[12pt,leqno,twoside]{article} 
\usepackage{bezier,amsmath,amssymb,stmaryrd,lscape}

\input xy
\xyoption{all}
\input{rk.sty}
\begin{document}

\title{A sufficient criterion for homotopy cartesianess}
\author{Alberto Canonaco, Matthias K\"unzer}
\maketitle

\begin{small}
\begin{quote}
\begin{center}{\bf Abstract}\end{center}\vspace*{2mm}
Suppose given a commutative quadrangle in a Verdier triangulated category such that there exists an induced isomorphism on the horizontally taken cones. Suppose that the endomorphism ring of the 
initial or the terminal corner object of this quadrangle satisfies a finiteness condition. Then this quadrangle is homotopy cartesian.
\end{quote}
\end{small}

\renewcommand{\thefootnote}{\fnsymbol{footnote}}
\footnotetext[0]{MSC 2000: 18E30.}
\renewcommand{\thefootnote}{\arabic{footnote}}

\begin{footnotesize}
\renewcommand{\baselinestretch}{0.7}%
\parskip0.0ex%
\tableofcontents%
\parskip1.2ex%
\renewcommand{\baselinestretch}{1.0}%
\end{footnotesize}%

\setcounter{section}{-1}

\section{Introduction}

In an abelian category, a commutative quadrangle is called {\it bicartesian} if its diagonal sequence is short exact, i.e.\ if it is a pullback and a pushout. A commutative quadrangle is 
bicartesian if and only if we get induced isomophisms on the horizontal kernels and on the horizontal cokernels.

In a triangulated category in the sense of {\sc Verdier} \bfcite{Ve63}{Def.\ 1-1}, a commutative quadrangle is called {\it homotopy cartesian} (or a {\it Mayer-Vietoris square}, or a 
{\it distinguished weak square}), if its diagonal sequence fits into a distinguished triangle. A homotopy cartesian square has a (non-uniquely) induced isomorphism
on the horizontally taken cones \bfcite{Ne01}{Lem.\ 1.4.4}. We consider the converse question\,: a commutative quadrangle that has an isomorphism induced on the horizontally taken cones, is it homotopy 
cartesian? We show this to be true if the endomorphism ring of the object in the terminal or initial corner satisfies a finiteness condition.

This finiteness condition is for instance satisfied for the endomorphism rings occurring in $\DD^\text{b}(A\modl)$, where $A$ is a finite-dimensional algebra over some field; or in $\ul{A\modl}$, 
where $A$ is a finite-dimensional Frobenius algebra over some field.

This finiteness condition, however, in general fails for the endomorphism rings occurring in $\KK^\text{b}(\Z\projl)$. We show by an example that the conclusion on our commutative 
quadrangle to be homotopy cartesian fails there as well. 
\section{A ring theoretical lemma}

Let $R$ be a ring. Denote by $\JJ(R)$ its Jacobson radical. 

If $R/\JJ(R)$ is artinian, we fix the notation $R/\JJ(R) \iso \prod_{i = 1}^n D_i^{k_i\ti k_i}$ for its Wedderburn decomposition, where $D_i$ is a skewfield for $1\le i\le n$.

A ring $R$ shall be called {\it head-finite} if its head $R/\JJ(R)$ is artinian and if in the Wedderburn decomposition of $R/\JJ(R)$, the skewfield $D_i$ is finite dimensional over its centre 
for each \mb{$1\le i\le n$}.

For example, finite dimensional algebras over some field are head-finite. For another example, a local ring $R$ for which $R/\JJ(R)$ is commutative is head-finite.

\begin{Lemma}
\label{Lem1}
Suppose given a head-finite ring $R$ and an element $\eps\in R$.

\begin{itemize}
\item[{\rm (1)}] There exists $\al\,\in\, R$ such that $1 + \eps + \al\eps^2$ is a unit in $R$.
\item[{\rm (2)}] There exists $\be\,\in\, R$ such that $1 + \eps + \eps^2\be$ is a unit in $R$.
\end{itemize}
\end{Lemma}

{\it Proof.} Assertion (2) follows by an application of (1) to $R^\0$, so it remains to prove (1).

Since an element $\rh\in R$ is a unit in $R$ if and only if 
$\rh + \JJ(R)$ is a unit in $R/\JJ(R)$, we may assume $\JJ(R) = 0$ and $R$ to be a product of matrix rings over skew fields which are finite-dimensional over their centres. 
Furthermore, we may assume $R$ to be a single matrix ring over a skewfield which is finite-dimensional over its centre $K$. In particular, we may assume $R$ to be a finite dimensional $K$-algebra.

Let $m\ge 0$ and $s(X)\in K[X]$ be such that $\eps$ is a root of the polynomial $X^m + X^{m+1}s(X)$. Let $\al := s(\eps)$. Then 
\[
(\eps + \al\eps^2)^{m+1} \= \big(\eps^m (1 +  s(\eps)\eps)\big)\big(\eps (1 + s(\eps)\eps)^m\big) \= 0\; ,
\]
and thus $1 + \eps + \al\eps^2$ is a unit in $R$.\qed

\begin{Remark}
\label{Rem1_5}\indent\rm
\begin{itemize}
\item[(1)] The conclusions of Lemma \ref{Lem1} do not hold for all rings, as the example $R = \Z$ and $\eps = 3$ shows.
\item[(2)] The conclusions of Lemma \ref{Lem1} hold for a local ring $R$, regardless whether $R/\JJ(R)$ is finite-dimensional over its centre or not.
\item[(3)] In Lemma \ref{Lem1}.(1), we do not claim that $\al\eps = \eps\al$. Whereas this property can be achieved if $R$ is a finite dimensional algebra over some field, as we have
seen in the proof of loc.\ cit., the first reduction step at the beginning of this proof possibly might not respect this property. 
\end{itemize}
\end{Remark}

\section{The criterion for homotopy cartesianess}

Let $\Cl$ be a triangulated category in the sense of {\sc Verdier} \bfcite{Ve63}{Def.\ 1-1}.

A commutative quadrangle 
\[
\xymatrix{
B\ar[d]_b\ar[r]^g & C\ar[d]^c \\
B'\ar[r]^{g'}     & C' \\
}
\]
in $\Cl$ is said to be {\it homotopy cartesian} if there exists a distinguished triangle containing the sequence
\[
\xymatrix@C=12mm{
B\ar[r]^(0.35){\smatze{b}{g}} & B'\ds C\ar[r]^(0.55){\smatez{g'}{-c}} & C' \; , \\
}
\]
cf.\ \bfcite{Ne01}{Def.\ 1.4.1}. We remark that by \bfcite{Ne01}{Lem.\ 1.4.4}, such a homotopy cartesian square fits into a morphism of
distinguished triangles of the form $(b,c,1)$ (and, by symmetry, also in one of the form $(g,g',1)$).

\begin{Proposition}
\label{Prop2}
Suppose given a commutative diagram in $\Cl$
\[
\xymatrix{
A\ar[r]^f\ar[d]^a_\wr & B\ar[r]^g\ar[d]^b & C\ar[r]^h\ar[d]^c & A[1]\ar[d]^{a[1]}_\wr \\
A'\ar[r]^{f'}         & B'\ar[r]^{g'}     & C'\ar[r]^{h'}     & A'[1]                 \\
}
\]
whose rows are distinguished triangles. 

\begin{itemize}
\item[{\rm (1)}] Suppose $\End_\Cl C'$ to be head-finite. Then the quadrangle $(g,g',b,c)$ is homotopy cartesian.
\item[{\rm (2)}] Suppose $\End_\Cl B $ to be head-finite. Then the quadrangle $(g,g',b,c)$ is homotopy cartesian.
\end{itemize}
\end{Proposition}

{\it Proof.} By duality, it suffices to prove (1). By isomorphic replacement at $A'$, we may assume that $A = A'$ and $a = 1$.

By \bfcite{Ne01}{Lem.\ 1.4.3}, there exists $C\lraa{\w c} C'$ such that the quadrangle $(g,g',b,\w c)$ is homotopy cartesian and such that $h'\0 \w c = h$. 
It suffices to show that the quadrangles $(g,g',b,c)$ and $(g,g',b,\w c)$ are isomorphic.

Since $(\w c - c)\0 g = 0$, there exists $A[1]\lraa{\psi} C'$ such that $\psi\0 h = \w c - c$. Let 
\[
\eps \; :=\; \psi\0 h' \;\in\; \End_\Cl C'\; . 
\]
By assumption on $\End_\Cl C'$, we may apply Lemma \ref{Lem1} to find an element $\al\in \End_\Cl C'$ such that $1 + \eps + \al\0\eps^2$ is a unit in $\End_\Cl C'$, i.e.\ an automorphism of $C'$.

We claim that we have the following isomorphism of commutative quadrangles.
\[
\xymatrix{
B\ar[rr]^g\ar[dr]^b\ar@{=}[dd] &                                & C\ar[dr]^c\ar@{=}'[d][dd] &                                          \\
                               & B'\ar[rr]^(0.4){g'}\ar@{=}[dd] &                           & C'\ar[dd]^{1+\eps + \al\,\0\,\eps^2}_\wr \\
B\ar'[r]^(0.7)g[rr]\ar[dr]^b   &                                & C\ar[dr]^{\w c}           &                                          \\
                               & B'\ar[rr]^{g'}                 &                           & C'                                       \\
}
\]
In fact, 
\[
\barcl
(1 + \eps + \al\0\eps^2)\0 c 
& = & c + \psi\0 h'\0 c + \al\0\psi\0 h'\0\psi\0 h'\0 c \\ 
& = & c + \psi\0 h + \al\0\psi\0 h'\0\psi\0 h \\
& = & \w c + \al\0\psi\0 h'\0 (\w c - c) \\
& = & \w c + \al\0\psi\0 (h - h) \\
& = & \w c\; , \\
\ea
\] 
and
\[
\barcl
(1 + \eps + \al\0\eps^2)\0 g' 
& = & g' + \psi\0 h'\0 g' + \al\0\psi\0 h'\0\psi\0 h'\0 g' \\ 
& = & g'\; . \\ 
\ea
\]
\qed

\begin{Question}[open]
\label{Qu2_5}
Suppose given a diagram as in Proposition {\rm\ref{Prop2}}. Is the cone of $b$ isomorphic to the cone of $c$? 
\end{Question}

Note that the isomorphism in question is not required to satisfy any commutativities. 
In \S\ref{SecCounter}, we give an example in which there is no isomorphism between these cones that is compatible with the diagram.

Note that the middle quadrangle of such a diagram is a weak square by the kernel-cokernel criterion applied in the Freyd category of $\Cl$; cf.\ e.g.\ \bfcite{Ku05}{\S A.6.3, Def.\ A.9, Lem.\ A.11}.

\begin{Question}[open]
\label{Qu2_7}
Suppose given a commutative quadrangle that has an isomorphism induced on the horizontally and on the vertically taken cones. Is it homotopy cartesian?
\end{Question}

Of course, provided the endomorphism ring of its object in the initial or terminal corner is head-finite, such a quadrangle is homotopy cartesian by Proposition \ref{Prop2}.
\section{A counterexample}
\label{SecCounter}

We shall give an example of a commutative quadrangle in a triangulated category that horizontally fits into a morphisms of triangles containing an isomorphism as in Proposition \ref{Prop2}, but that 
is not homotopy cartesian; somewhat worse still, vertically, it does not fit into a morphism of triangles containing an isomorphism. This will show that the head-finiteness conditions in 
Proposition \ref{Prop2} cannot be entirely dropped.

Let $\Cl := \KK^\text{b}(\Z\projl)$. As to sign conventions, the standard distinguished triangle on a morphism $X\lraa{f} Y$ in $\Cl$ is given as follows.
\[
\xymatrix{
\vdots\ar[d]                  & \vdots\ar[d]                                   & \vdots\ar[d]                                                                                        & \vdots\ar[d]               \\
X^i\ar[r]^{f^i}\ar[d]^{\de^i} & Y^i\ar[r]^(0.4){\smatze{1}{0}}\ar[d]^{\dell^i} & Y^i\ds X^{i+1}\ar[r]^(0.6){\smatez{0}{1}}\ar[d]^(0.45){\rsmatzz{\dell^i}{f^{i+1}}{0\;}{-\de^{i+1}}} & X^{i+1}\ar[d]^{-\de^{i+1}} \\
X^{i+1}\ar[r]^{f^{i+1}}\ar[d] & Y^{i+1}\ar[r]^(0.4){\smatze{1}{0}}\ar[d]       & Y^{i+1}\ds X^{i+2}\ar[r]^(0.6){\smatez{0}{1}}\ar[d]                                                 & X^{i+2}\ar[d]              \\
\vdots                        & \vdots                                         & \vdots                                                                                              & \vdots                     \\
}
\]
We shall allow ourselves to omit zero object entries when displaying complexes. By $\Z^{\ds m}$ we denote the direct sum of $m$ copies of $\Z$, where $m\ge 2$.

\begin{Lemma}
\label{Lem2}
The following triangles are distinguished in $\Cl$ for every $a,b\in\Z$.

\[
\ba{l}
\xymatrix@C=24mm{
0\ar[r]\ar[d]             & 0\ar[r]\ar[d]                  & \Z\ar[r]^{1}\ar[d]^{\smatze{b}{a^2}}       & \Z\ar[d]^{a^2} \\
\Z\ar[d]^{-a^2}\ar[r]^{b} & \Z\ar[d]\ar[r]^{\smatze{1}{0}} & \Z^{\ds 2}\ar[d]\ar[r]^{\smatez{0\;}{\;1}} & \Z\ar[d]       \\
\Z\ar[r]                  & 0\ar[r]                        & 0\ar[r]                                    & 0              \\
}\\
\ea
\leqno (1)
\]

\[
\ba{l}
\xymatrix@C=24mm{
0\ar[r]\ar[d]                         & \Z\ar[r]^1\ar[d]^{\rsmatze{-a}{a^2}}            & \Z\ar[r]^a\ar[d]^{\rsmatze{-a^3}{a^2}}      & \Z\ar[d]^{a^2} \\
\Z\ar[d]^{-a^2}\ar[r]^{\smatze{1}{0}} & \Z^{\ds 2}\ar[d]\ar[r]^{\smatzz{a^2}{0}{0}{1}}  & \Z^{\ds 2}\ar[d]\ar[r]^{\smatez{-1\;}{\;0}} & \Z\ar[d]       \\
\Z\ar[r]                              & 0\ar[r]                                         & 0\ar[r]                                     & 0              \\
} \\
\ea
\leqno (2)
\]

\[
\ba{l}
\xymatrix@C=24mm{
0\ar[r]\ar[d]              & \Z\ar[r]^{\smatze{1}{a}}\ar[d]^{\rsmatze{-a^3}{a^2}} & \Z^{\ds 2}\ar[r]^{\smatez{a\;}{\;-1}}\ar[d]^{\smatzz{a^2}{0}{0}{a^2}} & \Z\ar[d] \\
\Z\ar[r]^{\smatze{a^2}{0}} & \Z^{\ds 2}\ar[r]^{\rsmatzz{0}{1}{-1}{0}}             & \Z^{\ds 2}\ar[r]                                                      & 0        \\
}\\
\ea
\leqno (3)
\]

\[
\ba{l}
\xymatrix@C=24mm{
0\ar[r]\ar[d]  & 0\ar[r]\ar[d] & \Z\ar[r]^{1+a}\ar[d]^{a^2} & \Z\ar[d] \\
\Z\ar[r]^{a^2} & \Z\ar[r]^{1-a}  & \Z\ar[r]                 & 0        \\
}\\
\ea
\leqno (4)
\]

\end{Lemma}

{\it Proof.} It is enough to show that each triangle is of the form
\[
\xymatrix@C=12mm{
X\ar[r]^f & Y\ar[r]^{u\0 g} & Z'\ar[r]^{h\0 u^{-1}} & X[1]
}
\]
with $X\lraa{f} Y\lraa{g} Z \lraa{h} X[1]$ a standard distinguished
triangle and $Z\lraa{u} Z'$ an isomorphism in $\Cl$. Indeed, (1) is already
standard, and it is straightforward to check that in cases (2), (3)
and (4) one can take respectively the following morphisms for $u$.
\[
\xymatrix@C=24mm{
\Z^{\ds 2}\ar[r]^{\smatez{1}{0}}\ar[d]_{\smatdz{-a}{1}{\; a^2}{0}{0}{a^2}} & \Z\ar[d]^{\rsmatze{-a^3}{a^2}} \\
\Z^{\ds 3}\ar[r]^{\rsmatzd{a^2\!\!}{\;\; 0}{-1}{0}{1}{0}}                  & \Z^{\ds 2}                     \\
}
\]
\[
\xymatrix@C=24mm{
\Z^{\ds 2}\ar[r]^{\rsmatzz{1}{0}{a}{-1}}\ar[d]_{\smatzz{-a^3}{a^2}{\phantom{-}a^2}{0}} & \Z^{\ds 2}\ar[d]^{\smatzz{a^2}{0}{0}{a^2}} \\
\Z^{\ds 2}\ar[r]^{\rsmatzz{0}{1}{-1}{\; 0}}                                            & \Z^{\ds 2}                                 \\
}
\]
\[
\xymatrix@C=24mm{
\Z\ar[r]^{1-a}\ar[d]_{a^2} & \Z\ar[d]^{a^2} \\
\Z\ar[r]^{1-a}             & \Z             \\
}
\]
\qed

Let $a\in\Z$ be such that $a\ge 3$. Consider the following morphism of distinguished triangles. The differentials of the complexes are displayed from lower left to upper right, and the 
triangles are displayed from left to right. Notice that the triangles
are distinguished because they can be obtained by applying axiom (TR 2)
to the triangles (2) and (1) (with $b=-a$) of Lemma \ref{Lem2}.

{\footnotesize
\[
\ba{l}
\xymatrix@!@C=5mm{
                             &                                                                     &\Z^{\ds 2}\ar[rrr]^{\smatzz{a^2}{0}{0}{1}}\ar@{=}'[d]'[dd][ddd] &                          &                                                                           &\Z^{\ds 2}\ar[rrr]^{\smatez{-1}{0}}\ar'[d]'[dd]^{\smatez{0}{1}}[ddd] &                          &                                                                                   &\Z\ar[rrr]\ar'[d]'[dd][ddd] &                                                           &                              & 0\ar@{=}[ddd] \\
                             &\Z\ar[rrr]^1\ar[ur]^*+<-2mm,-2mm>{\smatze{-a}{a^2}}\ar@{=}'[d][ddd]  &                                                                &                          &\Z\ar[rrr]^a\ar[ur]^*+<-2mm,-2mm>{\rsmatze{-a^3}{a^2}}\ar'[d][ddd]^(0.25)1 &                                                                     &                          &\Z\ar[rrr]^{\rsmatze{-1}{0}}\ar[ur]^*+<-1mm,0mm>{\scm a^2}\ar'[d][ddd]^(0.25){1+a} &                            &                                                           &\Z^{\ds 2}\ar[ur]\ar@{=}[ddd] &               \\
0\ar[rrr]\ar[ur]\ar@{=}[ddd] &                                                                     &                                                                & 0\ar[rrr]\ar[ur]\ar[ddd] &                                                                           &                                                                     & 0\ar[rrr]\ar[ur]\ar[ddd] &                                                                                   &                            &\Z\ar[ur]^(0.4)*+<-2mm,-2mm>{\smatze{a}{-a^2}}\ar@{=}[ddd] &                              &               \\
                             &                                                                     &\Z^{\ds 2}\ar'[r]'[rr]^{\smatez{0}{1}}[rrr]                     &                          &                                                                           &\Z\ar'[r]'[rr][rrr]                                                  &                          &                                                                                   & 0\ar'[r]'[rr][rrr]         &                                                           &                              & 0             \\
                             &\Z\ar'[rr]^(0.75)1[rrr]\ar[ur]_(0.55)*+<-2mm,-2mm>{\smatze{-a}{a^2}} &                                                                &                          &\Z\ar'[rr]^(0.75)a[rrr]\ar[ur]^*+<-1mm,0mm>{\scm a^2}                      &                                                                     &                          &\Z\ar'[rr]^(0.7){\rsmatze{-1}{0}}[rrr]\ar[ur]                                      &                            &                                                           &\Z^{\ds 2}\ar[ur]             &               \\
0\ar[rrr]\ar[ur]             &                                                                     &                                                                & 0\ar[rrr]\ar[ur]         &                                                                           &                                                                     & 0\ar[rrr]\ar[ur]         &                                                                                   &                            &\Z\ar[ur]_(0.55)*+<-2mm,-2mm>{\smatze{a}{-a^2}}            &                              &               \\
} \\
\ea
\leqno (\ast)
\]
}

As an aside, we remark that the second morphism from the left in $(\ast)$ is split epimorphic.

We claim that the middle quadrangle of $(\ast)$ is not homotopy cartesian. 

We {\it assume} the contrary. The diagonal sequence of the middle quadrangle is given as follows.
\[
\xymatrix@C=24mm{
\Z\ar[r]^{\smatze{1}{a}}\ar[d]_{\rsmatze{-a^3}{a^2}} & \Z^{\ds 2}\ar[r]^{\smatez{a\;}{\;-1-a}}\ar[d]^{\smatzz{a^2}{0}{0}{a^2}} & \Z\ar[d] \\
\Z^{\ds 2}\ar[r]^{\rsmatzz{0}{1}{-1}{0}}             & \Z^{\ds 2}\ar[r]                                                        & 0        \\
}
\]
In contrast, the following sequence fits into the distinguished
triangle (3) of Lemma \ref{Lem2}.
\[
\xymatrix@C=24mm{
\Z\ar[r]^{\smatze{1}{a}}\ar[d]_{\rsmatze{-a^3}{a^2}} & \Z^{\ds 2}\ar[r]^{\smatez{a\;}{\;-1}}\ar[d]^{\smatzz{a^2}{0}{0}{a^2}} & \Z\ar[d] \\
\Z^{\ds 2}\ar[r]^{\rsmatzz{0}{1}{-1}{0}}             & \Z^{\ds 2}\ar[r]                                                      & 0        \\
}
\]
By uniqueness of the cone up to isomorphism in $\Cl$, there exists a commutative triangle in $\Cl$
\[
\xymatrix{
                                                                                                         & & &                                      & \Z\ar[d] \\
\Z^{\ds 2}\ar[d]_{\smatzz{a^2}{0}{0}{a^2}}\ar[rrr]_{\smatez{a\;}{\;-1-a}}\ar[urrrr]^{\smatez{a\;}{\;-1}} & & & \Z\ar[ur]_*+<-2mm,0mm>{\scm s}\ar[d] & 0        \\
\Z^{\ds 2}\ar[rrr]\ar[urrrr]|(0.73)\hole                                                                 & & & 0\ar[ur]                             &          \\
}
\]
with $s\in\{-1,+1\}$. We conclude that $sa\con_{a^2} a$ and $-s - sa \con_{a^2} -1$. If $s = 1$, then the second congruence gives $a \con_{a^2} 0$, which is impossible since $a\ge 2$.
If $s = -1$, then the first congruence gives $2a\con_{a^2} 0$, which is impossible since $a\ge 3$. We have arrived at a {\it contradiction}.

We claim that vertically, the middle quadrangle of $(\ast)$ does not
fit into a morphism of distinguished triangles that contains an
isomorphism. We note that this claim implies the preceding claim
(cf.\ \mb{\bfcite{Ne01}{Lem.\ 1.4.4}}), which we will thus have proven twice.

We {\it assume} the contrary. Inserting triangles (1) (with $b=-a^3$)
and (4) of Lemma \ref{Lem2} vertically, we obtain the following commutative quadrangle above our given one, where $t\in\{-1,+1\}$. Differentials are displayed from lower left to 
upper right.
\[
\xymatrix@C=12mm{
                                                           &\Z\ar[rr]^t\ar'[d]^(0.6)*+<-1mm,0mm>{\smatze{1}{0}}[dd] &                                  &\Z\ar[dd]^{1-a} \\
0\ar[ur]\ar[rr]\ar[dd]                                     &                                                        & 0\ar[ur]\ar[dd]                  &                \\
                                                           &\Z^{\ds 2}\ar[rr]^(0.35){\smatez{-1}{0}}                &                                  &\Z              \\
\Z\ar[ur]^(0.6)*+<-1mm,-1mm>{\rsmatze{-a^3}{a^2}}\ar[rr]^a &                                                        &\Z\ar[ur]_*+<-1mm,-1mm>{\scm a^2} &                \\
}
\]
Commutativity of this quadrangle in $\Cl$ means that $t(1-a) \con_{a^2} -1$. If $t = -1$, then $a \con_{a^2} 0$ ensues, which is impossible since $a\ge 2$. If $t = 1$, then $2 \con_{a^2} a$ and
hence $2\con_a 0$ ensues, which is impossible since $a\ge 3$. We have arrived at a {\it contradiction}.

\parskip0.0ex
\begin{footnotesize}

\parskip1.2ex

\vspace*{3cm}

\begin{minipage}{8cm}
Alberto Canonaco\\
Universit\`a degli studi di Pavia\\
Dipartimento di matematica F.\ Casorati\\
Via Ferrata, 1\\
I-27100 Pavia\\
alberto.canonaco@unipv.it\\
\end{minipage}
\hfill
\begin{minipage}{8cm}
\begin{flushright}
Matthias K\"unzer\\
Lehrstuhl D f\"ur Mathematik\\
RWTH Aachen\\
Templergraben 64\\
D-52062 Aachen \\
kuenzer@math.rwth-aachen.de \\
www.math.rwth-aachen.de/$\sim$kuenzer\\
\end{flushright}
\end{minipage}
\end{footnotesize}
\end{document}